\documentclass[12pt]{article}

\usepackage[reqno]{amsmath}
\usepackage{amssymb, amsthm}
\usepackage{enumerate,arydshln/arydshln, url} 
\newtheoremstyle{plainsl}%
	{\topsep}
	{\topsep}
	{\slshape} 
	{}
	{\normalfont\bfseries}
	{.}
	{ }
	{}

\swapnumbers 

\theoremstyle{plainsl}
\newtheorem{theorem}{Theorem}[section]
\newtheorem{lemma}[theorem]{Lemma}
\newtheorem{corollary}[theorem]{Corollary}
\newtheorem{prop}[theorem]{Proposition}

\newcommand\one{{\bf1}}
\DeclareMathOperator{\spanof}{span}

\DeclareMathOperator\sym{Sym}
\DeclareMathOperator\aut{Aut}

\newcommand\re{{\mathbb R}}

\renewcommand\proof{\noindent\textsl{Proof. }}
\newcommand\sqr[2]{{\vbox{\hrule height.#2pt
    \hbox{\vrule width.#2pt height#1pt \kern#1pt
        \vrule width.#2pt}\hrule height.#2pt}}}

\renewcommand\qed{%
	\ifmmode\eqno\sqr53
	\else\nolinebreak\ \hfill\sqr53\medbreak\fi}

\numberwithin{equation}{section}

\newcommand{\remove}[1]{}

\newcommand\ind[2]{{\mathrm{ind}_{#2}({#1})} }
\newcommand\res[2]{{\mathrm{res}_{#2}({#1})} }
\newcommand\indg[2]{{\mathrm{ind}_{#2}(1_{#1})} }

\title{An algebraic proof of the Erd\H{o}s-Ko-Rado theorem for intersecting families of perfect matchings}

\author{Chris Godsil \footnote{Research supported by NSERC.}\\
\small  Department of Combinatorics and Optimization \\[-0.8ex]
\small University of Waterloo,  Waterloo, Ontario, Canada\\[-0.8ex]
\small \texttt{cgodsil@math.uwaterloo.ca}\\
Karen Meagher \small{*}\\ [-0.8ex]
\small  Department of Mathematics and Statistics \\[-0.8ex]
\small University of Regina,  Regina, Saskatchewan, Canada\\[-0.8ex]
\small \texttt{karen.meagher@uregina.ca}
}

\begin{document}
\maketitle

\abstract{In this paper we give a proof that the largest set of
  perfect matchings, in which any two contain a common edge, is the
  set of all perfect matchings that contain a fixed edge. This is a
  version of the famous Erd\H{o}s-Ko-Rado theorem for perfect
  matchings. The proof given in this paper is algebraic, we first
  determine the least eigenvalue of the perfect matching derangement
  graph and use properties of the perfect matching polytope. We also
  prove that the perfect matching derangement graph is not a Cayley
  graph.}

\section{Introduction}

A perfect matching in the complete graph $K_{2k}$ is a set of $k$
vertex disjoint edges.  Two perfect matchings intersect if they
contain a common edge. In this paper we use an algebraic method to prove
that the natural version of the Erd\H{o}s-Ko-Rado (EKR) theorem holds for
perfect matchings. This theorem shows that the largest set of perfect
matchings, with the property that any two intersect, is the set of the
all perfect matchings that contain a specific edge.

The algebraic method given in this paper is similar to the proof given
in~\cite{MR2489272} that the natural version of the EKR theorem holds
for permutations. In this paper we determine the least eigenvalue for
the perfect matching derangement graph. This, with the
Delsarte-Hoffman bound, implies that a maximum intersecting set of
perfect matchings corresponds to a facet in the perfect matching
polytope. The characterization of the maximum set of intersecting
perfect matchings follows from the characterization of the facets of
this polytope. Further, with this characterization, we are able to
prove that the perfect matching derangement graph is not a Cayley graph.

Meagher and Moura~\cite{MR2156694} proved a version of the EKR theorem
holds for intersecting uniform partitions using a counting
argument~\cite{MR2156694}. This result includes the EKR theorem for
perfect matchings. It is interesting that the counting argument
in~\cite{MR2156694} is straight-forward, except for the case of
perfect matchings; in this case a more difficult form of the counting
method is necessary.

\section{Perfect matchings}

A \textsl{perfect matching} is a set of vertex disjoint edges in the
complete graph $K_{2k}$. This is equivalent to a partition of a
set of size $2k$ into $k$-disjoint classes, each of size $2$. The number of
perfect matchings in $K_{2k}$ is
\[
\frac{1}{k!} \binom{2k}{2} \binom{2k-2}{2} \cdots \binom{2}{2} =
(2k-1)(2k-3) \cdots 1.
\]
For an odd integer $n$ define
\[
n!! = n(n-2)(n-4) \cdots 1,
\]
thus, there are $(2k-1)!!$ perfect matchings.

We say that two perfect matchings are \textsl{intersecting} if they both
contain a common edge, and a set of perfect matchings is intersecting if
the perfect matchings in the set are pairwise intersecting. If $e$
represents a pair from $\{1,\dots,2k\}$, then $e$ is an edge of
$K_{2k}$. Define $S_e$ to be the set of all perfect matchings that
include the edge $e$.  Clearly the sets $S_e$ are
intersecting. The sets $S_e$ are called the \textsl{canonically
  intersecting} sets of perfect matchings. For every $e$
\[
|S_e| = (2k-3)!! = (2k-3)(2k-5) \cdots 1.
\]

The main result of this paper can be stated as follows.

\begin{theorem}\label{thm:matching}
  The largest set of intersecting perfect matchings in $K_{2k}$ has
  size $(2k-3)!!$. The only sets that meet this bound are the
  canonically intersecting sets of perfect matchings.
\end{theorem}

\section{Perfect matching derangement graph}

One approach to proving EKR theorems for different objects is
to define a graph where the vertices are the objects and two objects
are adjacent if and only if they are not intersecting (see
\cite{MR2489272, MR2739502, MR2231096} for
just a few examples of where this is done). This is the
approach that we take with the perfect matchings.

We use the standard graph notation. A \textsl{clique} in a graph is a
set of vertices in which any two are adjacent, a \textsl{coclique} is
a set of vertices in which no two are adjacent. If $X$ is a graph,
then $\omega(X)$ denotes the size of the largest clique, and
$\alpha(X)$ is the size of the largest coclique. A graph is
\textsl{vertex transitive} if its automorphism group is transitive
on the vertices. In this case, there is a relationship between
the maximum clique size and maximum coclique size known as the
\textsl{clique-coclique bound}.

\begin{theorem}\label{thm:cliquecoclique}
Let $X$ be a vertex-transitive graph, then
\[
\alpha(X) \, \omega(X) \leq |V(X)|.\qed
\]
\end{theorem}

The \textsl{eigenvalues} of a graph are the eigenvalues of the
adjacency matrix of the graph. Similarly, the eigenvectors and
eigenspaces of the graph are the eigenvectors and eigenspaces of the
adjacency matrix.

Define the \textsl{perfect matching derangement graph} $M(2k)$ to be the graph
whose vertices are all perfect matchings on $K_{2k}$ and vertices are
adjacent if and only if they have no edges in common.
Theorem~\ref{thm:matching} is equivalent to the statement that the
size of the maximum coclique in $M(2k)$ is $(2k-3)!!$ and that only the
canonically intersecting sets, $S_e$, meet this bound.

The number of vertices in $M(2k)$ is $(2k-1)!!$.  The degree of
$M(2k)$, denoted by $d(2k)$, is the number of perfect matchings that
do not contain any the edges from some fixed perfect matching. This
number can be calculated using the principle of inclusion-exclusion:
\begin{align} \label{eq:matchingdegree}
d(2k) = \sum_{i=0}^{k-1}(-1)^{i} \binom{k}{i}(2k\!-\!2i\!-\!1)!!.
\end{align}
In practise, this formula can be tricky to use, but we will make use the
following simple lower bound on $d(2k)$.

\begin{lemma}\label{lem:matchinglowerbound}
For any $k$
\[
d(2k) >  (2k-1)!! - \binom{k}{1}(2k-3)!! .
\]
\end{lemma}
\proof
For any $i \in \{0, \dots, k-1\}$
\begin{align*}
\binom{k}{i}(2k-2i-1)!! &= \frac{i+1}{k-i} \binom{k}{i+1}(2k-2i-1) (2k-2(i+1)-1)!! \\
  &>\binom{k}{i+1}(2k-2(i+1)-1)!!
\end{align*}
(since $\frac{i+1}{k-i} (2k-2i-1) >1$ for these values of $i$).
This implies that the terms in Equation~\ref{eq:matchingdegree} are
strictly decreasing in absolute value. Since
it is an alternating sequence, the first two terms give a lower bound
on $d(2k)$.
\qed

Next we give some simple properties of the perfect matching
derangement graph, including a simple proof of the bound in
Theorem~\ref{thm:matching} that uses the clique-coclique bound.

\begin{theorem}\label{thm:properties}
Let $M(2k)$ be the perfect matching derangement graph.
\begin{enumerate}
\item The graph $M(2k)$ is vertex transitive, and $\sym(2k)$ is a subgroup of the
  automorphism group of $M(2k)$.
\item The size of a maximum clique in $M(2k)$ is $2k-1$.
\item The size of a maximum coclique in $M(2k)$ is
$(2k-3)!!$.
\end{enumerate}
\end{theorem}
\proof It is clear that the group $\sym(2k)$ acts transitively on the
perfect matchings and, though this action, each permutation in
$\sym(2k)$ gives an automorphism of $M(2k)$. 

Let $C$ be a clique in $M(2k)$. For every perfect matching in $C$, the element
$1$ is matched with a different element of $\{2,3,\dots, 2k\}$.  Thus
the size of $C$ is no more than $2k-1$. A $1$-factorization of the
complete graph on $2k$ vertices is a clique of size
$\frac{1}{k}\binom{2k}{2}$ in $M(2k)$. Since a $1$-factorization of $K_{2k}$
exists for every $k$, the size of the maximum clique is exactly
$\frac{1}{k}{2k \choose 2} = 2k-1$.

Since $M(2k)$ is vertex transitive, the clique-coclique bound,
Theorem~\ref{thm:cliquecoclique}, holds so
\begin{align*}
\alpha(M(2k)) & \leq \frac{(2k-1)!!}{\frac{1}{k}{2k \choose 2}}  = (2k-3)!!.
\end{align*}
Since the set $S_e$ meets this bound $\alpha(M(2k)) = (2k-3)!!$.
\qed

\section{Perfect matching association scheme}

We have noted that the group $\sym(2k)$ acts on the set of perfect matchings. Under
this action, the stabilizer of a single perfect matching is isomorphic
to the wreath product of $\sym(2)$ and $\sym(k)$. This is a subgroup
of $\sym(2k)$ and is denoted by $\sym(2) \wr \sym(k)$. Thus the set of perfect
matchings in $K_{2k}$ correspond to the set of cosets 
\[
\sym(2k)/(\sym(2) \wr \sym(k)).
\]
This implies that the action of $\sym(2k)$ on the perfect matchings is
equivalent to the action of $\sym(2k)$ on the cosets
$\sym(2k)/(\sym(2) \wr \sym(k))$.  This action produces a permutation
representation of $\sym(2k)$. We will not give much detail on the
representation theory of the symmetric group, rather we will simply
state the results that we need and refer the reader to any standard text
on the representation theory of the symmetric group, such
as~\cite{MR1153249} or \cite{MR1093239}.

Each irreducible representation of $\sym(2k)$ corresponds to an
integer partition $\lambda \vdash 2k$; these representations will be
written as $\chi_\lambda$. Information about the representation is
contained in the partition. For example, the dimension of the
representation can be found just from the partition using the
\textsl{hook length formula}.

For any group $G$, the trivial representation of $G$ is denoted by
$1_G$.  If $\chi$ is a representation of a group $H \leq \sym(n)$,
then $\ind{\chi}{\sym(n)}$ is the representation of $\sym(n)$ induced
by $\chi$. Similarly, if $\chi$ is a representation of $\sym(n)$, then
$\res{\chi}{H}$ is the restriction of $\chi$ to $H$. The permutation
representation of $\sym(n)$ acting on $\sym(2k)/(\sym(2) \wr \sym(k))$
is the representation induced on $\sym(2k)$ by the trivial
representation on $\sym(2) \wr \sym(k)$. This representation is
denoted by $\indg{\sym(2) \wr \sym(k)}{\sym(2k)}$ (see \cite[Chapter
13]{EKRbook} for more details).

For an integer partition $\lambda \vdash k$ with
$\lambda=(\lambda_1,\lambda_2, \dots, \lambda_\ell)$, let $2\lambda$
denote the partition $(2\lambda_1,2\lambda_2, \dots, 2\lambda_\ell)$
of $2k$.  It is well-known (see, for example, \cite[Example
2.2]{MR627512}) that the decomposition of the permutation representation
of $\sym(2k)$ from its action on the perfect matchings is
\[
\indg{\sym(2) \wr \sym(k)}{\sym(2k)} = \sum_{\lambda \vdash k} \chi_{2\lambda}.
\]

The multiplicity of each irreducible representation in this
decomposition is one, this implies that $\indg{\sym(2) \wr
  \sym(k)}{\sym(2k)}$ is a \textsl{multiplicity-free representation}.
This implies that the adjacency matrices of the orbitals from the
action of $\sym(2k)$ on the cosets $\sym(2k) / (\sym(2) \wr \sym(k))$
defines an association scheme on the perfect matchings (see
\cite[Section 13.4]{EKRbook} for more details and a proof of this
result). This association scheme is known as the \textsl{perfect
  matching scheme}. Each class in this scheme is labelled with a
partition $2\lambda = (2\lambda_1,2\lambda_2, \dots,
2\lambda_\ell)$. Two perfect matchings are adjacent in a class if
their union forms a set of $\ell$ cycles with lengths $2\lambda_1,
2\lambda_2, \dots, 2\lambda_\ell$ (this association scheme is described
in more detail in~\cite[Section 15.4]{EKRbook} and \cite{MR1285207}).

The graph $M(2k)$ is the union of all the classes in this association
scheme in which the corresponding partition contains no part of size
two. This means that each eigenspace of $M(2k)$ is the union of
$\chi_{2\lambda}$-modules of $\sym(2k)$ where $\lambda \vdash k$. If
$\xi$ is an eigenvalue of $M(2k)$, and its eigenspace includes the
$\chi_{2\lambda}$-module, then we say that $\xi$ is the eigenvalue
belonging to the $\chi_{2\lambda}$-module. Conversely, we denote the
eigenvalue belonging to the $\chi_{2\lambda}$-module by
$\xi_{2\lambda}$.

The next lemma contains a formula to calculate the eigenvalue
belonging to the $\chi_{2\lambda}$-module. This gives considerable
information about the eigenvalues of $M(2k)$. For a proof the general
form of this formula see \cite[Section 13.8]{EKRbook}, we only state
the version specific to perfect matchings. If $M$ denotes a perfect
matching and $\sigma \in \sym(2k)$, we will use $M^\sigma$ to denote
the matching formed by the action of $\sigma$ on $M$.

\begin{lemma}\label{lem:evalueorbital}
  Let $M$ be a fixed perfect matching in $K_{2k}$.  Let $H \subseteq
  \sym(2k)$ be the set of all elements $\sigma \in \sym(2k)$ such that
  $M$ and $M^\sigma$ are not intersecting.  The eigenvalue of $M(2k)$
  belonging to the $\chi_{2\lambda}$-module is
\[
\xi_{2\lambda}   =\frac{d(2k)}{2^k \, k!} \sum_{x \in  H}\chi_{2\lambda}(x). 
\]
The $\chi_{2\lambda}$-module is a subspace of the
$\xi_{2\lambda}$-eigenspace and the dimension of this subspace is
$\chi_{2\lambda}(1)$. \qed
\end{lemma}

This formula can be used to calculate the eigenvalue corresponding to
a module for the matching derangement graph. This formula is not
effective to determine all the eigenvalues for a general matching
derangement graph. In Section~\ref{sec:evaluesMatching} we will show
another way to find some of the eigenvalues.

\section{Delsarte-Hoffman bound}
\label{sec:DHBound}

In Section~\ref{sec:evaluesMatching}, we will give an alternate proof
of the bound in Theorem~\ref{thm:matching} that uses the eigenvalues
of the matching derangement graph. This proof is based on the
Delsarte-Hoffman bound, which is also known as the ratio
  bound. The advantage of this bound is that when equality holds we
get additional information about the cocliques of maximum size. This
information can be used to characterize all the sets that meet the bound. The
Delsarte-Hoffman bound is well-known and there are many references, we
offer~\cite[Theorem 2.4.1]{EKRbook} for a proof.

\begin{theorem}\label{thm:ratcoq}
  Let $X$ be a $k$-regular graph with $v$ vertices and let $\tau$ be
  the least eigenvalue of $A(X)$.  Then
\[
  \alpha(X) \le\frac{v}{1-\frac{k}{\tau}}.
\]
If equality holds for some coclique $S$ with characteristic vector
$v_S$, then
\[
v_S -\frac{|S|}{|V(X)|}\one
\]
is an eigenvector with eigenvalue $\tau$.\qed
\end{theorem}

If equality holds in the Delsarte-Hoffman bound, we say that the
maximum cocliques are \textsl{ratio tight}.  

The Delsarte-Hoffman bound can be used to prove the EKR theorem for
sets. Similar to the situation for the perfect matchings, the group
$\sym(n)$ acts on the subsets of  $\{1,\dots,n\}$ of size $k$.
This action is equivalent to the action of $\sym(n)$ on the cosets $\sym(n)
/(\sym(n-k)\times \sym(k))$. This action corresponds to a permutation
representation, namely
\begin{align}
\indg{\sym(n-k) \times \sym(k)}{\sym(n)} = \sum_{i=0}^{k} \chi_{[n-i,i]}. \label{eq:decomp}
\end{align}
(Details can be found in any standard text on the representation
theory of the symmetric group.)  This representation is multiplicity
free and the orbital schemes from this action
is an association scheme better known as the \textsl{Johnson scheme}. 

The \textsl{Kneser graph} $K(n,k)$ is the graph whose
vertices are all the $k$-sets from $\{1,\dots, n\}$ and two vertices are
adjacent if and only if they are disjoint. The Kneser graph is a graph in the
Johnson scheme, it is the graph that corresponds to the orbitals of
pairs of sets that do not intersect.  A coclique in $K(n,k)$ is a set
of intersecting $k$-sets. The Kneser graph is very well-studied and
all of its eigenvalues are known (see~\cite[Chapter 7]{MR1829620} or
\cite[Section 6.6]{EKRbook} for a proof).

\begin{prop}\label{prop:evaluesKneser}
The eigenvalues of $K(n,k)$ are
\[
(-1)^{i}\binom{n-k-i}{r-i}
\]
with multiplicities $\binom{n}{i} -\binom{n}{i-1}$
for $i\in \{0,\dots, k\}$.
\end{prop}

If we apply the Delsarte-Hoffman bound to $K(n,k)$ we get the following theorem
which is equivalent to the standard EKR theorem. The characterization
follows from the second statement in the Delsarte-Hoffman bound,
see~\cite[Section 6.6]{EKRbook} for details.

\begin{theorem}
  Assume that $n >2k$. The size of the largest coclique in $K(n,k)$ is
 $ \binom{n-1}{k-1}$, and the only cocliques of this size consist of
  all $k$-sets that contain a common fixed element.\qed
\end{theorem}

To apply the Delsarte-Hoffman bound to $M(2k)$, we first need to determine
the value of the least eigenvalue of $M(2k)$. We do not calculate all
the eigenvalues of $M(2k)$, rather we calculate the two eigenvalues
with the largest absolute value and then show that all other eigenvalues
have smaller absolute value.

\section{Dimensions of the representations of $\sym(2k)$}

In this section we determine a list of the irreducible representations of
$\sym(2k)$ with small degree.  We will need to use the
\textsl{branching rule}; this is a rule that describes which
irreducible representations are in $\res{\chi_\lambda}{\sym(n-1)}$ and $\ind{\chi_\lambda}{\sym(n+1)}$
based on the structure of $\lambda$. For a proof of
this rule see~\cite[Corollary 3.3.11]{MR2643487}.

\begin{lemma}\label{lem:branching}
Let $\lambda \vdash n$, then
\[
\res{\chi_\lambda}{\sym(n-1)} = \sum \chi_{\lambda^-},
\]
where the sum is taken over all partitions $\lambda^{-}$ of $n-1$ that have a
Young diagram which can be obtained
by the deletion of a single box from the Young diagram of $\lambda$.
Further,
\[
\ind{\chi_\lambda}{\sym(n+1)} = \sum \chi_{\lambda^+},
\]
where the sum is taken over partitions $\lambda^+$ of $n+1$ that have a
Young diagram which can be obtained by the addition of a single box to
Young diagram of $\lambda$. \qed
\end{lemma}

\begin{lemma}\label{lem:eightspecial}
  For $n \geq 9$, let $\chi$ be a representation of $\sym(n)$ with
  degree less than $(n^2-n)/2$. If $\chi_\lambda$ is a constituent of
  $\chi$, then $\lambda$ is one of the following partitions of $n$:
\[
[n],\, [1^n], \, [n-1,1],\, [2,1^{n-2}],\, [n-2,2], \, [2,2,1^{n-4}],\, [n-2,1,1], \, [3,1^{n-3}].
\]
\end{lemma}
\proof We prove this by induction. The result for $n=9,10$ can be read
directly from the character table for $\sym(n)$ (these character
tables are available in~\cite{MR644144}, or from the GAP character
table library~\cite{GAP4}). We assume that the lemma holds for $n$ and
$n-1$.

Assume that $\chi$ is a representation of $\sym(n+1)$ that has
dimension less than 
\[
\frac{(n+1)^2 -(n+1)}{2} = (n^2 + n)/2, 
\]
but does not have one of the eight
irreducibles representations listed in the statement of the theorem as
a constituent.

Consider the restriction of $\chi$ to $\sym(n)$, this representation
will be denoted by $\res{\chi}{\sym(n)}$. If one of the eight
irreducible representations of $\sym(n)$ with dimension less than
$(n^2-n)/2$ is a constituent of $\res{\chi}{\sym(n)}$, then, using
the branching rule, we can determine the constituents of $\chi$. The
first column of Table~\ref{tab:tableofsmallreps} gives the irreducible
representations of $\sym(n)$ that have dimension less than
$(n^2-n)/2$. The second column lists the irreducible representations
that must be constituents of $\chi$, if the representation in the
first column is a constituent of $\res{\chi}{\sym(n)}$. These are
determined by the branching rule.

\renewcommand{\arraystretch}{1.25}
\begin{figure}[h!]
\begin{center}
\begin{tabular}{|l|l|} \hline
Constituent of $\res{\chi}{\sym(n)}$  & Constituents of $\chi$ \\ \hline
$[n]$ & $[n+1]$, $[n,1]$\\ 
$[n-1,1]$ & $[n,1]$, $[n-1,2]$, $[n-1,1,1]$ \\
$[n-2,2]$ & $[n-1,2]$, $[n-2,3]$, $[n-2,2,1]$ \\
$[n-2,1,1]$ & $[n-1,1,1]$, $[n-2,2,1]$, $[n-2,1,1,1]$ \\
$[1^n]$ & $[2,1^{n-1}]$, $[1^{n+1}]$ \\
$[2,1^{n-2}]$ & $[3,1^{n-2}]$, $[2,2,1^{n-3}]$, $[2,1^{n-1}]$ \\
$[2,2, 1^{n-4}]$ & $[3,2,1^{n-4}]$, $[2,2,2,1^{n-5}]$, $[2,2,1^{n-3}]$ \\
$[3, 1^{n-3}]$ & $[4,1^{n-3}]$, $[3,2,1^{n-4}]$, $[3,1^{n-2}]$ \\ \hline
\end{tabular}
\caption{Constituents of $\chi$, if $\res{\chi}{\sym(n)}$ has a constituent with degree less than $(n^2-n)/2$.}\label{tab:tableofsmallreps}
\end{center}
\end{figure}
\renewcommand{\arraystretch}{1}

Using Table~\ref{tab:tableofsmallreps}, we see that either one of the
eight representations of $\sym(n+1)$ in the statement of the theorem
is a constituent of $\chi$; or, using the hook length formula, the
dimension of $\chi$ is larger than $(n^2+n)/2$ (provided
that $n>8$, see Table~\ref{tab:mediumreps}). Either case is a
contradiction. So we can assume that $\res{\chi}{\sym(n)}$ does not
contain any of the representations from $\sym(n)$ with dimension less
than $(n^2-n)/2$.

\renewcommand{\arraystretch}{1.25}
\begin{figure}[h!]
\begin{center}
\begin{tabular}{|l|l|} \hline
Representation & Degree \\ \hline
$[n-2,3] $& $n(n+1)(n-4)/6$ \\
$[n-2,2,1] $& $(n+1)(n-1)(n-3)/3$\\
$[n-2,1,1,1] $& $n(n-1)(n-2)/6$\\
$[3,2,1^{n-4}] $& $(n+1)(n-1)(n-3)/3$ \\
$[2,2,2,1^{n-5}] $&  $n(n+1)(n-4)/6$ \\
$[4,1^{n-3}] $& $n(n-1)(n-2)/6$ \\ \hline
\end{tabular}
\caption{Degrees of the representations from Table~\ref{tab:tableofsmallreps} that
are larger than $(n^2+n)/2$.}\label{tab:mediumreps}
\end{center}
\end{figure}
\renewcommand{\arraystretch}{1}

If the decomposition of $\res{\chi}{\sym(n)}$ contains two irreducible
representations of $\sym(n)$, neither of which is one of the eight
irreducible representations with dimension less $(n^2-n)/2$, then the
dimension of $\chi$ must be at least
\[
2((n^2-n)/2) =n^2-n. 
\]
But since $n>3$, this is strictly larger than $(n^2+n)/2$. 

Thus $\res{\chi}{\sym(n)}$ must be an irreducible representation of
$\sym(n)$. The branching rules then implies that if $\chi_\lambda$ is a
constituent of $\chi$, then there is only one block whose removal from
the Young diagram of $\lambda$ will produce a Young diagram for a partition
on $n$. (The irreducible representation of $\sym(n)$ corresponding to
this Young diagram is equal to $\res{\chi}{\sym(n)}$.) The only such
Young diagrams are rectangular, so the constituents of $\chi$
must be of the form $\chi_{[s^t]}$, for some $s$ and $t$. Finally,
since $n>2$, there can only be one constituent of $\chi$ (or else
$\res{\chi}{\sym(n)}$ will not be irreducible), thus
we can assume that $\chi = \chi_{[s^t]}$.

Next consider the restriction of $\chi = \chi_{[s^t]}$ to $\sym(n-1)$,
this representation will be denoted by $\res{\chi}{\sym(n-1)}$. By the
branching rule, this can contain only the irreducible representations
of $n-1$ that correspond to the partitions $\lambda' = [s^{t-1}, s-2]$ and
$\lambda'' = [s^{t-2}, s-1, s-1]$. 

If $\lambda' = [n-1]$, then $s=1$ or $s=2$, and $s = n-1$; but this
contradicts the assumption that $n \geq 10$. If $\lambda' = [n-3,2]$,
then $t=2$ and $s=4$; if $\lambda' = [2,2,1^{n-5}]$, then $t=4$ and
$s=2$. Again, these cases cannot happen, since $n = st$ and $n$ is
assumed to be greater than $10$.  If $\lambda'$ is any of the other
eight partitions of $n-1$ that correspond to irreducible
representations of $\sym(n-1)$ dimension less that
$((n-1)^2-(n-1))/2$, then $s=3$ and $t = 1$; which again is a
contradiction with $n=st$ and $n >10$. Similarly, $\lambda''$ cannot
be any of the partitions corresponding to the eight representations of
$\sym(n-1)$ that have dimension less than $\left((n-1)^2
  -(n-1)\right)/2$.

Thus the dimension of $\chi$ is at least
\[
2 \frac{ (n-1)^2 -(n-1) }{2} = (n-1)^2 -(n-1)
\]
which is strictly greater than $(n^2 + n )/2$ for any $n \geq
7$.  \qed

\section{Eigenvalues of the matching derangement graph}
\label{sec:evaluesMatching}

In this section we determine the largest and the least eigenvalue of the matching
derangement graph and identify the modules to which they belong. First
we will use a simple method to show that these two values are
eigenvalues of $M(2k)$.

For any edge $e$ in $K_{2k}$, the partition $\pi = \{S_e, V(M(2k)) /
S_e\}$ is an equitable partition of the vertices in $M(2k)$. In fact,
$\pi$ is the orbit partition formed by the stabilizer of the edge $e$
in $\sym(2k)$ (this subgroup is isomorphic to $\sym(2) \times
\sym(2k-2)$) acting on the set of all vertices of $M(2k)$. The
quotient graph of $M(2k)$ with respect to this partition is
\[
M(2k)/\pi = 
    \begin{pmatrix}
        0 & d(2k) \\
       \frac{1}{2k-2}\, d(2k) & \frac{2k-3}{2k-2}\, d(2k)
   \end{pmatrix}.
\]
The eigenvalues
for the quotient graph $M(2k)/\pi$ are 
\[
	d(2k),\qquad -\frac{d(2k)}{2k-2}.
\]
Since $\pi$ is equitable, these are also eigenvalues of $M(2k)$. The
next result identifies which modules these eigenvalues belong to.

\begin{lemma}\label{lem:2matchingevalues}
  The eigenvalue of $M(2k)$ belonging to the $\chi_{[2k]}$-module is $d(2k)$, and the
  eigenvalue of $M(2k)$ belonging to the $\chi_{[2k-2,2]}$-module is $-d(2k)/(2k-2)$.
\end{lemma}
\proof The first statement is clear using the formula in
Lemma~\ref{lem:evalueorbital}.

To prove the second statement, we will consider the equitable
partition $\pi$ defined above. The partition $\pi$ is
the orbit partition of $\sym(2) \times \sym(2k-2)$ acting on the
perfect matchings. Let $H =\sym(2) \times \sym(2k-2)$ and denote the
cosets of $H$ in $\sym(2k)$ by $\{x_0H = H, x_1H, \dots, x_{2k^2-k-1}H\}$.

The $-d(2k)/(2k-2)$-eigenvector of $M(2k)/\pi$ lifts to an eigenvector
$v$ of $M(2k)$.  A simple calculation shows that the entries of $v$
are $1-\frac{1}{2k-1}$ or $-\frac{1}{2k-1}$, depending on if the index
of the entry is in $S_e$, or not. This means that $v =v_e
-\frac{1}{2k-1} \one$, where $v_e$ is the characteristic vector of
$S_e$.

The group $\sym(2k)$ acts on the edges of $K_{2k}$, and for each
$\sigma \in \sym(2k)$, we can define
\[
v^\sigma = v_{e^\sigma} -\frac{1}{2k-1} \one.
\]
The vector $v$ is fixed by any permutation in $H$ under this action. 
If we define
\[
V =\spanof\{v^\sigma \,:\, \sigma \in \sym(n) \},
\]
then $V$ is a subspace of the $-d(2k)/(2k-2)$-eigenspace. Moreover $V$
is invariant under the action of $\sym(2k)$, so it is also a
$\sym(2k)$-module. To prove this lemma we need to show that $V$ is
isomorphic to the $\chi_{[2k-2,2]}$-module.

Let $W$ be the $\sym(2k)$-module for the induced representation
$\indg{H}{\sym(2k)}$.  By Equation~\ref{eq:decomp}, $W$ is the
sum of irreducible modules of $\sym(2k)$ that are isomorphic to
$M_{[2k]}$, $M_{[2k-1,1]}$ and $M_{[2k-2,2]}$. The vector space $W$ is
isomorphic to the vector space of functions $f \in L(\sym(n))$ that
are constant on $H$. For each coset $xH$, the function
$\delta_{xH}(\sigma)$ is equal to $1$ if $\sigma$ is in $xH$ and $0$
otherwise. The functions $\delta_{xH}$ form a basis for $W$.

Define the map $f$ so that
\[
f(v^\sigma) = \delta_{\sigma H} - \frac{1}{2k-1} \sum_{i =0}^{2k^2-k-1} \delta_{x_i H}. 
\]
Since $v^\sigma = v^\pi$ if and only if $\sigma H = \pi H$, this
function is well-defined. Further, it is a $\sym(2k)$-module
homomorphism. Thus $V$ is isomorphic to a submodule of $W$. Since $V$
is not trivial, it must be the $\chi_{[2k-2,2]}$-module, since it is the only
module (other than the trivial) that is common to both
$\indg{H}{\sym(2k)}$ and $\indg{\sym(2) \wr \sym(k)}{\sym(2k)}$.\qed

Next we will bound the size of the other eigenvalues.  This bound
follows from the straightforward fact that if $A$ is the adjacency
matrix of a graph, then the trace of the square of $A$ is equal to
both the sum of the squares of the eigenvalues of $A$, and to twice
the number of edges in the graph. The proof of this result closely
follows the proof the least eigenvalue of the derangement graph of the
symmetric group by Ellis~\cite{MR2876315}.

\begin{theorem}\label{thm:leastevaluematching}
  For $\lambda \vdash k$, the absolute value of the eigenvalue of $M(2k)$
  corresponding $\chi_{2\lambda}$-module is strictly less than
  $d(2k)/(2k-2)$, unless $\lambda = [2k]$ or $\lambda =[2k-2,2]$.
\end{theorem}
\proof Let $A$ be the adjacency matrix of $M(2k)$ and use
$\xi_{2\lambda}$ to denote the eigenvalue for the $\chi_{2\lambda}$-module.

The sum of the eigenvalues of $A^2$ is twice the number of edges in $M(2k)$, that is
\[
\sum_{\lambda \vdash k}\chi_{2\lambda}(1) \xi^2_{2\lambda} = (2k-1)!! \, d(2k).
\]
From Lemma~\ref{lem:2matchingevalues} we know the eigenvalues for two of the
modules, so this bound can be expressed as 
\[
\sum_{\stackrel{\lambda \vdash k}{\lambda \neq [k], [k-1,1]}}\chi_{2\lambda}(1) \xi^2_{2\lambda}  = 
(2k\!-\!1)!! \, d(2k) - d(2k)^2 - ( 2k^2\!-\!3k) \left(\frac{d(2k)}{2k\!-\!2}\right)^2.
\]
Since all the terms in left-hand side of the above summation are
positive, any single term is less than the sum. Thus 
\[
\chi_{2\lambda}(1) \xi^2_{2\lambda}  \leq  (2k\!-\!1)!! \, d(2k) - d(2k)^2 - ( 2k^2\!-\!3k) \left(\frac{d(2k)}{2k\!-\!2}\right)^2,
\]
(where $\lambda \vdash k$ and $\lambda \neq [k], [k-1,1]$).
If $|\xi_{2\lambda}|  \geq d(2k)/(2k-2)$, then this reduces to
\[
\chi_{\lambda}(1) \leq  \frac{ (2k-1)!! (2k-2)^2}{d(2k)} - 6k^2 +11k-4.
\]
Using the bound in Lemma~\ref{lem:matchinglowerbound}, this implies that
\[
\chi_{\lambda}(1) <  2k^2-k = \frac{(2k)^2-(2k)}{2}.
\]

If $|\xi_{2\lambda}| \leq d(2k)/(2k-2)$, then $2\lambda$ must be one
of the eight irreducible representations in Lemma~\ref{lem:eightspecial}. Thus
$2\lambda$ must be one of $[2k]$ and $[2k-2,2]$, which
proves the result.  \qed

We restate this result in terms of the least eigenvalue of the
matching derangement graph; noting that
Theorem~\ref{thm:leastevaluematching} implies that only the
$\chi_{[2k-2,2]}$-module has $-d(2k)/(2k-2)$ as its eigenvalue.

\begin{corollary}\label{cor:minevaluematching}
	The smallest eigenvalue of $M(2k)$ is $-d(2k)/(2k-2)$ and
        the multiplicity of this eigenvalue is $2k^2-3k$. \qed
\end{corollary}

\section{ $\chi_{[2k-2,2]}$-module}
\label{sec:module}

Applying the Delsarte-Hoffman bound with the fact that
$-d(2k)/(2k-2)$ is the least eigenvalue of $M(2k)$, proves that
the cocliques $S_e$ are ratio tight since
\[
	\frac{|V(M(2k)|}{1-\frac{d}{\tau}} =
	\frac{(2k-1)!!}{1-\frac{d(2k)}{-\frac{d(2k)}{2k-2}}} = (2k-3)!!.
\]

For $S$ a maximum coclique in $M(2k)$ we will use $v_S$ to denote the
characteristic vector of $S$. The ratio bound implies that $|S| =
(2k-3)!!$ and further that
\[
v_s - \frac{1}{2k-1}\one
\]
is a $-d(2k)/(2k-2)$-eigenvector. This vector is called the
\textsl{balanced} characteristic of $S$, since is it orthogonal to
that all ones vector. Since the $\chi_{[2k-2,2]}$-module is the only module
for which the corresponding eigenvalue is the least (this follows
directly from Theorem~\ref{thm:leastevaluematching}) we have the
following result which will be used to determine the structure of the
maximum cocliques in $M(2k)$.

\begin{lemma}\label{lem:inmodule}
  The characteristic vector for any maximum coclique in $M(2k)$ is in
  the direct sum of the $\chi_{[2k]}$-module and the
  $\chi_{[2k-2,2]}$-module. \qed
\end{lemma}

A perfect matching is a subset of the edges in the complete
graph, and thus can be represented as a characteristic vector; this is
a vector in $\re^{\binom{2k}{2}}$. Define the \textsl{incidence
  matrix} for the perfect matchings in $K_{2k}$ to be the matrix $U$
whose rows are the characteristic vectors of the perfect matchings of
$K_{2k}$. The columns of $U$ are indexed by the edges in the complete
graph and the rows are indexed by the perfect matchings. The column of
$U$ corresponding to the edge $e$ is the characteristic vector of the
canonical intersecting set of matchings $S_e$.

We will show that the characteristic vector of any maximum coclique of
$M(2k)$ is a linear combination of the columns of $U$.

\begin{lemma}\label{lem:matchingbalancedspan} 
The characteristic vectors of the canonical cocliques of $M(2k)$ span
the direct sum of the  $\chi_{[2k]}$-module and the $\chi_{[2k-2,2]}$-module. 
\end{lemma}
\proof Let $v_e$ be the characteristic vector of $S_e$. From
Lemma~\ref{lem:inmodule}, the vector $v_e - \frac{1}{2k-1}\one$ is in the
$\chi_{[2k-2,2]}$-module, and $v_e$ is in the direct sum of the
$\chi_{[2k]}$-module and the $\chi_{[2k-2,2]}$-module. So all that
needs to be shown is that the span of all the vectors $v_e$ has dimension
$2k^2-3k+1$, or equivalently, that the rank of $U$ is $2k^2-3k+1$.

Let $I$ denote the $\binom{k}{2} \times \binom{k}{2}$ identity matrix
and $A(2k,2)$ the adjacency matrix of the Kneser graph $K(2k,2)$. Then
\[
U^TU = (2k-3)!! I + (2k-5)!! A(2k,2).
\]
By Proposition~\ref{prop:evaluesKneser}, $0$ is an eigenvalue of this
matrix with multiplicity $2k -1$. Thus the rank of $U^TU$ (and hence
$U$) is $\binom{2k}{2} -(2k-1) = 2k^2-3k+1$.
 \qed

 Putting this result with the comments at the beginning of this
 section, we have the following corollary.

\begin{corollary}\label{cor:linearcombo}
The characteristic vector of a maximum coclique in the perfect
matching derangement graph is in the column space of $U$. \qed 
\end{corollary}

Next we will show that this implies that any maximum coclique is a
canonical coclique. To do this we will consider a polytope based on
the perfect matchings.

\section{The perfect matching polytope}

The convex hull of the set of characteristic vectors for all the
perfect matchings of a graph $K_{2k}$ is called the \textsl{perfect
  matching polytope} of $K_{2k}$.  Let $U$ be the incidence matrix
defined in the previous section, then the perfect matching polytope is
the convex hull of the rows of $U$.  A \textsl{face} of the perfect
matching polytope is the convex hull of the rows where $Uh$ achieves
its maximum for some vector $h$. A \textsl{facet} is a maximal proper
face of a polytope.

If $S$ is a maximum coclique in $M(2k)$, then from
Corollary~\ref{cor:linearcombo}, we know that $Uh = v_s$ for some
vector $h$. If a vertex of $K_{2k}$ is in $S$, then the corresponding
row of $Uh$ is equal to $1$; conversely, if a vertex of $K_{2k}$ is
not in $S$, then the corresponding row of $Uh$ is equal to $0$. Thus a
maximum intersecting set of perfect matchings is a facet of the
perfect matching polytope. In this section, we will give a
characterization of the facets of the perfect matching polytope for
the complete graph.

Let $S$ be a subset of the vertices of $K_{2k}$ and define the
\textsl{boundary} of $S$ to be the set of edges that join a vertex in
$S$ to a vertex not in $S$. The boundary is denoted by $\partial S$
and is also known as an \textsl{edge cut}.  If $S$ is a subset of the
vertices of $K_{2k}$ of odd size, then any perfect matching in $K_{2k}$ must
contain at least one edge from $\partial S$. If $S$ is a
single vertex, then any perfect matching contains exactly one element of
$\partial S$. It is an amazing classical result of Edmonds that these
two constraints characterize the perfect matching polytope for any
graph. For a proof of this result see Schrijver~\cite{Schrijver1983}.

\begin{theorem}\label{thm:schrijver}
	Let $X$ be a graph. A vector $x$ in $\re^{|E(X)|}$ lies in the perfect
	matching polytope of $X$ if and only if:
	\begin{enumerate}[(a)]
	\item $x\ge0$;\label{eq:schrijvera}
        \item if $S=\{u\}$ for some $u\in V(X)$, then 
              $\sum_{e \in \partial S}x(e)=1$; \label{eq:schrijverb}
	\item if $S$ is an odd subset of $V(X)$ with $|S| \geq 3$, then 
	      $\sum_{e\in\partial S} x(e)\ge1$.\label{eq:schrijverc}
	\end{enumerate}
	If $X$ is bipartite, then $x$ lies in the perfect matching
    polytope if and only if the first two conditions hold.\qed
\end{theorem}

The constraints in Equation~(\ref{eq:schrijverb}) define an affine
subspace of $\re^{|E(X)|}$. The perfect matching polytope is the
intersection of this subspace with affine half-spaces defined by the
conditions in Equation~(\ref{eq:schrijvera}) and
Equation~(\ref{eq:schrijverc}); hence the points in a proper face of the
polytope must satisfy at least one of these conditions with equality.

For any graph $X$ (that is not bipartite) the vertices of a facet
are either the perfect matchings that miss a given edge, or the
perfect matchings that contain exactly one edge from $\partial S$ for
some odd subset $S$.

It follows from Theorem~\ref{thm:schrijver} that every perfect
matching in $K_{2k}$ is a vertex in the perfect matching polytope for
the complete graph. But we can also determine the vertices of every
facet in this polytope.

\begin{lemma}\label{lem:odd-facet}
	In the matching polytope of $K_{2k}$, the vertices of a facet of
	maximum size are the perfect matchings that do not contain a given
	edge.
\end{lemma}
\proof Let $F$ be a facet of the polytope of maximum size. From the
above comments, equality holds in at least one of equations
\[
\sum_{e\in\partial S} x(e)\ge1
\]
for all $x \in F$.  Suppose $S$ is the subset that defines such an
equation, then $S$ is an odd subset of the vertices in $K_{2k}$ for
which $\sum_{e\in\partial S} x(e) = 1$ for all $x \in F$.

Let $s$ be the size of $S$. Each perfect matching with exactly
one edge in $\partial S$ consists of the following: a matching of size $(s-1)/2$
covering all but one vertex of $S$; an edge joining this missed vertex
of $S$ to a vertex in $\overline{S}$; and a matching of size $(2k-s-1)/2$
covering all but one vertex in $\overline{S}$.  
Hence there are
\[
	(s-2)!! \; s(2k-s) \; (2k-s-2)!! = s!!(2k-s)!!
\]
such perfect matchings. We denote this number by $N(s)$ and observe that
\[
	\frac{N(s-2)}{N(s)} = \frac{2k-s+2}{s}.
\]
Hence for all $s$ such that $3\le s \le k$ we see that the values
$N(s)$ are strictly decreasing, so the maximum size of a set of such
vertices is $N(3)=3(2k-3)!!$.

On the other hand, the number of perfect matchings in $K_{2k}$ that do
not contain a given edge is
\[
	(2k-1)!!-(2k-3)!! = (2k-2) (2k-3)!!.
\]
Since this is always larger than $N(3)$, the lemma follows.\qed

We now have all the tools to show that any maximum intersecting set of
perfect matchings is the set of all matchings that contain a fixed edge.

\begin{theorem}\label{lem:facetEKRcomplete}
  The largest coclique in $M(2k)$ has size $(2k-3)!!$. The only cocliques
  that meet this bound are the canonically intersecting sets of
  perfect matchings.
\end{theorem}
\proof Let $S$ be a maximum coclique in $M(2k)$ and let $v_S$ be the
characteristic vector of $S$. Then $|S| = (2k-3)!!$, by the
Delsarte-Hoffman bound and Corollary~\ref{cor:minevaluematching}. The
Delsarte-Hoffman bound, along with
Theorem~\ref{thm:leastevaluematching}, further imply that the vector
$v_S - \frac{1}{2k-1}\one$ is in the $\chi_{[2k-2,2]}$-module.

By Lemma~\ref{lem:matchingbalancedspan}, $v_S - \frac{1}{2k-1}\one$ is
a linear combination of the balanced characteristic vectors of the
canonical cocliques. This also implies that $v_S$ is a linear
combination of the characteristic vectors of the canonical cocliques.
So there exists a vector $x$ such that $Ux =v_S$ (where $U$ is the
matrix defined in Section~\ref{sec:module}).

Finally, by Lemma~\ref{lem:odd-facet}, $\overline{S}$ is a face of maximal
size an it consists of all the perfect matching that avoid a fixed
edge.  This implies that $S$ is a canonical coclique of $M(2k)$.\qed

\section{$M(2k)$ is not a Cayley graph}

The \textsl{derangement graph} for a permutation group is the graph
with the elements of the group as its vertices, and with permutations
$\sigma$ and $\tau$ adjacent if $\tau\sigma^{-1}$ is a
derangement. (If the group is the symmetric group, then this graph is
usually just called the \textsl{derangement graph}.) For any group,
the derangement graph is a Cayley graph for the group with the set of
derangements as the connection set. The matching derangement graph can
be viewed as an analogue of the derangement graph.  But, in contrast, we
will show that it follows from Theorem~\ref{thm:matching} that $M(2k)$
is not a Cayley graph when $k\ge3$.

\begin{lemma}\label{lem:autgroup}
	If $k\ge3$ then $\aut(M(2k))\cong\sym(2k)$.
\end{lemma}

\proof 
From Theorem~\ref{thm:matching}, there are exactly $\binom{2k}{2}$
maximum cocliques in $M(2k)$---these are the canonical cocliques, and
each canonical coclique is the set of all perfect matchings that contain a
fixed edge from $K_{2k}$. If $\alpha\in\aut(M(2k))$, then it determines a permutation of the
$\binom{2k}2$ canonical cocliques, and thus determines a permutation
$\hat\alpha$ of the $\binom{2k}2$ edges of $K_{2k}$.

Since $\alpha$ is an automorphism of $M(2k)$, if $S_{\{a,b\}}$ and $S_{\{c,d\}}$ are 
canonical cliques, then
\[
	S_{\{a,b\}}^\alpha \cap S_{\{c,d\}}^\alpha = (S_{\{a,b\}}\cap S_{\{c,d\}})^\alpha.
\]
If $|\{a,b\} \cap \{c,d\}|=1$, then $S_{\{a,b\}}\cap S_{\{c,d\}} =
\emptyset$ and the above equation shows that $ S_{\{a,b\}}^\alpha \cap
S_{\{c,d\}}^\alpha = \emptyset$. This implies that
$|\{a,b\}^{\hat\alpha} \cap \{c,d\}^{\hat\alpha}|=1$. Similarly, if the
edges $\{a,b\}$ and $\{c,d\}$ are disjoint, then their images under $\hat\alpha$ are also
disjoint. 

Consider $L(K_{2k})$, the line graph of the graph $K_{2k}$ (this is
the graph with the edges of $K_{2k}$ as its vertices, and two vertices
are adjacent if the edges share a vertex of $K_{2k}$). What we have
shown is that $\hat\alpha$ is an automorphism of $L(K_{2k})$.  The
fact that $\aut(L(K_{2k}))\cong\aut(K_{2k})\cong\sym(2k)$ 
(see, for example, \cite[Theorem 5.3]{MR0130186}) 
completes the theorem.\qed

\begin{theorem}
	If $k\ge3$, then $M(2k)$ is not a Cayley graph.
\end{theorem}

\proof Assume that $M(2k)$ is a Cayley graph on the group $G$.  Then,
by Lemma~\ref{lem:autgroup}, $G$ must be a subgroup of $\sym(2k)$
acting regularly on the vertices of $M(2k)$. Since the number of
vertices in $M(2k)$ is odd, this implies that $|G|$ is also odd and
therefore, applying a big hammer, the group $G$ is solvable. If $p$
and $q$ are distinct primes and $p^a$ and $q^b$ are the largest powers
of $p$ and $q$ that divide $|G|$, then by Hall's theorem for solvable
groups (see~\cite[Theorem 9.3.1]{MR0103215}), we see that $G$ has a
subgroup of order $p^aq^b$.

We will show provided that $k \geq 3$, then it is
possible to choose distinct primes $p$ and $q$ so that
\[
    k \leq p,q < 2k.
\]
This implies that $G$ has a subgroup $H$ of order $pq$.  Since $p$ and
$q$ are both greater than $k$ and less than $2k$, then $p$ does not
divide $q-1$, and $q$ does not divide $p-1$. This implies that the
Sylow $p$- and $q$-subgroups of $H$ are normal, and so it follows that
$H$ is abelian.  Since $|H|$ is square-free, $H$ is cyclic
(see~\cite[Section 4.4]{MR2286236}).

Thus, if $M(2k)$ is a Cayley graph and we can find two such primes
$p,q$, then we can conclude that there is a cyclic subgroup of
$\sym(2k)$ with order $pq$.  But, any element of $\sym(2k)$ of order
$p$ or $q$ is a cycle of length, respectively, $p$ or $q$, we deduce
that $\sym(2k)$ does not contain a cyclic subgroup of order $pq$ and
that we have a contradiction.

The final problem is to show that it is possible to choose $p$ and $q$
as needed. In 1952, Jitsuro Nagura~\cite{MR-JN} proved that if $m \ge
25$, there is always a prime between $m$ and $(1 + 1/5)m$.  So if $k$
is at least $25$, there is a prime between $k$ and $6/5 k$, and
another prime between $6/5 k +1$ and $6/5(6/5k+1) < 2k$. For $k\geq
25$, we will use these two prime numbers for $p$ and $q$. If $k<25$,
the following table lists primes that can be used for $p$ and $q$.

\begin{figure}[h!]
\begin{center}
\begin{tabular}{ |r | c |} \hline
	$k$ & primes\\ \hline
	\hline
	18--24& 29, 31\\
	12--17& 19, 23\\
	7--11&  11, 13\\
	6&   7, 11\\
	4,5&   5,7\\
	3&     3,5	\\ \hline
\end{tabular}
\caption{Two distinct primes between $k$ and $2k$.}
\end{center}
\end{figure}

This completes the proof.\qed

\end{document}